\input amstex

\documentstyle{amsppt}

\refstyle{A}

\nologo

\hoffset .25 true in
\voffset .2 true in

\hsize=6.1 true in
\vsize=8.5 true in

\define\Pdot{\bold P^\bullet}
\define\Qdot{\bold Q^\bullet}
\define\Fdot{\bold F^\bullet}
\define\mf{F_{{}_{f, \bold p}}}
\define\mfz{F_{{}_{f, \bold 0}}}
\define\mfpn{F_{{}_{f, \bold p_0}}}

\topmatter

\abstract We show that there are obstructions to the existence of certain types of invariant subspaces of the Milnor monodromy; this
places restrictions on the cohomology of Milnor fibres of non-isolated hypersurface singularities.
\endabstract

\title Invariant Subspaces of the Monodromy \endtitle

\author David B. Massey \endauthor

\address{David B. Massey, Dept. of Mathematics, Northeastern University, Boston, MA, 02115, USA} \endaddress

\email{DMASSEY\@NEU.edu}\endemail



\endtopmatter

\document

\vskip .2in

\noindent\S0. {\bf Introduction}  

\vskip .2in

Let $\Cal U$ be an open neighborhood of the origin in $\Bbb C^{n+1}$, where $n\geqslant 2$.
Let
$f:(\Cal U,
\bold 0)\rightarrow(\Bbb C, 0)$ be a complex analytic function which has a smooth $1$-dimensional critical locus, $\Sigma f$, at the
origin. Note that this, combined with the fact that $n\geqslant 2$, implies that $f$ is reduced at the origin.

For $\bold p\in V(f)$, we are interested in the cohomology of the Milnor fibre $\mf$. For technical reasons, it is easier to work with
field coefficients, rather than with integer coefficients. Of course, by varying the base field through finite fields, one may still
detect torsion in the integral cohomology. Hence, throughout this paper, we fix a field $\frak K$, and all cohomology spaces will be
$\frak K$-vector spaces.

Since we care only about the analytic-type of the germ of $f$ at the origin, we may make an analytic change of coordinates and assume that 
$\Sigma f$ is, in fact, a portion of the complex line
$\Bbb C\times \{\bold 0\}$. Thus, we shall assume that $S:=\Sigma f = (\Bbb C\times\{\bold 0\})\cap\Cal U$. In addition,  we shall assume
that
$\Cal U$ has been chosen small enough so that
$S\subseteq V(f)$, and so that
$\big\{\Cal U-V(f), V(f)-S, S-\{\bold 0\}, \{\bold 0\}\big\}$ is a Whitney stratification. Finally, let us use $z_0$ for the first
coordinate of $\Bbb C^{n+1}$, and assume that
$\Cal U$ is small enough and that $z_0$ is generic enough so that $\Sigma(f_{|_{V(z_0)}})=\{\bold 0\}$.

\vskip .2in

At each point $\bold p\in V(f)-S$, the Milnor fibre $\mf$ is contractible; hence, at such points, the reduced  cohomology
$\widetilde H^*(\mf) = 0$. 

At each point $\bold p\in S-\{\bold 0\}$, the Milnor fibre is the cross-product of a disk and the Milnor fibre of a generic hyperplane
slice. Hence, at such points, 
$\widetilde H^*(\mf) = 0$, except in degree $n-1$; we shall denote the dimension of $\widetilde H^{n-1}(\mf)$   by
$\lambda^1_f$. (This is the $1$-dimensional L\^e number; see [{\bf M2}].) Finally,
$\widetilde H^*(\mfz)$ is possibly non-zero only in degrees
$n-1$ and $n$.

\vskip .2in

However, there are further restrictions on the cohomologies of the various $\mf$.  The cohomology groups $\widetilde H^{n-1}(\mf)$
actually form a local system over the punctured disk $S-\{\bold 0\}$; this local system is characterized by fixing a point $\bold p_0\in
S-\{\bold 0\}$ and considering the monodromy isomorphism $$\nu:
\widetilde H^{n-1}(\mfpn)\rightarrow \widetilde H^{n-1}(\mfpn),$$ which is obtained by letting $\bold p_0$ travel around the origin in
$S$.  We refer to this as the {\it internal monodromy} (it is sometimes referred to as the {\it vertical monodromy}). Note that this is
most definitely {\bf not} the Milnor monodromy, which results from letting the value of
$f$ travel around the origin. 

Even at the origin, the cohomology of the Milnor fibre is not independent of the local system determined by $\nu$. The relationships
between the cohomology of $\mfz$ and the cohomologies of the other $\mf$ are encoded by the complex of sheaves of {\it vanishing cycles},
$\phi_f\Bbb {\frak K}^\bullet_{{}_{\Cal U}}$ (see [{\bf D}], [{\bf K-S}]); this complex of sheaves is defined on all of $V(f)$, but
is supported on
$S$. It is more convenient to restrict the vanishing cycles to their support $S$ and, for technical reasons, to shift the complex into
non-positive degrees; hence, we let $\Pdot:=\big(\phi_f{\frak K}^\bullet_{{}_{\Cal U}}\big)_{|_S}[n]$.

Can we obtain, up to isomorphism, every bounded constructible complex of $\frak K$-vector spaces on a complex line by varying the function
$f$ through all analytic functions with line singularities and looking at the corresponding $\Pdot$'s? Certainly not. The discussion above
indicates restrictions on the degrees in which one may have non-zero cohomology; after our shift, we find that $\Pdot_{|_{S-\{\bold 0\}}}$
is a local system, concentrated in degree
$-1$ and characterized by $\nu$, while the stalk $\Pdot_\bold 0$ has possibly non-zero cohomology only in degrees $-1$ and $0$. All of
these properties are a reflection of the well-known fact that $\Pdot$ must be a {\bf perverse sheaf} on $S$.

In fact, the properties above almost completely describe perverse sheaves on $S$, which are constructible with respect to
$\{S-\{\bold 0\}, \{\bold 0\}\}$. The only additional required property is implied by the cosupport condition at the origin: the inclusion
map $i:S-\{\bold 0\}\hookrightarrow S$ must induce an injection
$$ H^{-1}(\Pdot)_\bold 0\cong \Bbb H^{-1}(S; \Pdot)\hookrightarrow \Bbb H^{-1}(S-\{\bold 0\}; \Pdot).
$$ Using that the hypercohomology of a local system (in degree $-1$) on a punctured disk is given by
$\Bbb H^{-1}(S-\{\bold 0\}; \Pdot)\cong\operatorname{ker}(\operatorname{id}-\nu)$ (and $\Bbb H^{0}(S-\{\bold 0\};
\Pdot)\cong\operatorname{coker}(\operatorname{id}-\nu)$), this extra condition can be stated as: $i$ induces an inclusion
$$
\widetilde H^{n-1}(\mfz)\hookrightarrow \operatorname{ker}(\operatorname{id}-\nu).
$$

\vskip .1in

The reasonable question to ask now is: can every {\bf perverse} sheaf on $S$ be obtained by varying the function $f$? Again, the answer is
``no''; the complex $\Pdot$ must be self-dual, i.e., $\Pdot$ must be isomorphic to its own Verdier dual; this is well-known (it follows
from non-natural isomorphisms $\Cal D\circ\phi_f[-1]\cong \phi_f[-1]\circ\Cal D$ and $\Cal D\big({\frak K}^\bullet_{{}_{\Cal
U}}[n+1]\big)\cong {\frak K}^\bullet_{{}_{\Cal U}}[n+1]$; see, for instance, [{\bf B}]). In Section 1, we shall use MacPherson and Vilonen's
results in [{\bf M-V}] to provide a description of the category of perverse sheaves on a line, and then we shall be able to explain
self-duality in a down-to-Earth manner. MacPherson and Vilonen's description is essential for the understanding the results of this paper.

\vskip .3in

The goal of this paper is to show that there are further general restrictions on the vanishing complex $\Pdot$, i.e., one can not even
obtain every self-dual perverse sheaf on $S$ by varying $f$; there are invariant subspaces of the monodromy isomorphism which impose
additional structure on
$\Pdot$.  This result implies new restrictions on the cohomology of the Milnor fibre of a line singularity in affine space.

\vskip .3in

However, despite the fact that our methods use the derived category and perverse sheaves, the statement of our main theorem is very
classical.

Given the complex analytic function $f:(\Cal U, \bold 0)\rightarrow(\Bbb C, 0)$, whose critical locus at the origin is a line, we may
select a generic linear form $L$, and consider the one-parameter family of isolated singularities given by $f_t:=f_{|_{V(L-L(t))}}$, and
we may consider the corresponding Milnor numbers $\mu_t$. For all small $t\neq 0$, the $\mu_t$ have a common value, namely $\lambda^1_f$.
The main theorem of this paper is Theorem 3.3:

\vskip .3in

\noindent{\bf Theorem}. {\it If $\mu_0=1+\lambda^1_f$, then there is a strict inequality $\operatorname{dim}_{\frak K}\widetilde
H^{n-1}(\mfz)<\lambda^1_f$. }

\vskip .3in

The previously-known general bound on the reduced cohomology in degree $n-1$ was the non-strict inequality $\operatorname{dim}_{\frak
K}\widetilde H^{n-1}(\mfz)\leqslant\lambda^1_f$.

\vskip .4in

\noindent\S1. {\bf The Category of Perverse Sheaves on a Line}  

\vskip .2in

We continue with the notation from the introduction.

\vskip .2in

In order to explain what self-duality implies, it will be useful to avail ourselves of the description of the category of perverse sheaves
on a line which is provided by MacPherson and Vilonen in [{\bf M-V}].

\vskip .1in

If $g:\Cal U\rightarrow\Bbb C$ is any complex analytic function, we may consider the functors $\psi_g[-1]$ and $\phi_g[-1]$ given by
shifting the nearby and vanishing cycles, respectively, i.e., for a complex
$\Fdot$ on $\Cal U$, $\psi_g[-1]\Fdot := (\psi_g\Fdot)[-1]$ and $\phi_g[-1]\Fdot := (\phi_g\Fdot)[-1]$. The reason for defining these
shifted functors is that $\psi_g[-1]$ and $\phi_g[-1]$ take perverse sheaves to perverse sheaves, i.e., they yield functors from the
Abelian category of perverse sheaves on $\Cal U$ to the Abelian category of perverse sheaves on $V(g)$.

\vskip .1in

Now, suppose  that $\Qdot$ is an arbitrary perverse sheaf on the line $S$, and let $\hat z_0:={z_0}_{|_{S}}$. Then, 
$\phi_{\hat z_0}[-1]\Qdot$ and $\psi_{\hat z_0}[-1]\Qdot$ are perverse sheaves on the single point $\bold 0$, and consequently have
possibly non-zero cohomology only in degree $0$. Let $W:= H^0\big(\phi_{\hat z_0}[-1]\Qdot\big)_\bold 0$ and $V:= H^0\big(\psi_{\hat
z_0}[-1]\Qdot\big)_\bold 0$. Then, there is the canonical map $\gamma:V\rightarrow W$ and the variation map
$\delta:W\rightarrow V$, and $\delta\circ\gamma=\operatorname{id}-\nu$, where $\nu$ is, again, the internal monodromy isomorphism of
$\Qdot_{|_{S-\{\bold 0\}}}$, which arises from considering what happens as the stalks travel around a small circle centered at the origin
in
$S$. 

Thus, a perverse sheaf on $L$ yields two vector spaces, $V$ and $W$, an isomorphism $\nu:V\rightarrow V$, and  maps $\gamma:V\rightarrow
W$ and 
$\delta:W\rightarrow V$ such that $\delta\circ\gamma=\operatorname{id}-\nu$. This situation is nicely represented by a commutative triangle

\vbox{
$$ V@>\hskip .1in\operatorname{id}-\nu\hskip .1in>>V
$$
\vskip  -.2in
$$
\hbox{}\hskip -0.0in\gamma\searrow\hskip .27in\nearrow\delta
$$
\vskip -.2in
$$
\hbox{}\hskip 1in W \hskip 1in .
$$ } We refer to such a commutative diagram as an {\it M-V triangle}. The category of perverse sheaves on $S$ (constructible with respect
to
$\{S-\{\bold 0\},
\{\bold 0\}\}$) is equivalent to the category of M-V triangles, where a morphism of M-V triangles is defined in the obvious way: a
morphism is determined by linear maps
$\tau:V\rightarrow V^\prime$ and
$\eta: W\rightarrow W^\prime$ such that 
$$ V\hskip .05in@>\ \gamma\ >> \hskip .05in W @>\ \hskip .05in\delta\ >> \hskip .05in V 
$$
$$
\hbox{}\hskip -.15in\tau\downarrow\hskip .34in\eta\downarrow\hskip .31in\tau\downarrow
$$
$$ V^\prime @>\ \gamma^\prime\ >> W^\prime @>\ \delta^\prime\ >> V^\prime
$$ commutes.

Note that one can recover the stalk cohomology of $\Qdot$ at the origin from the associated M-V triangle:
$$ H^{-1}(\Qdot)_\bold 0\cong \operatorname{ker}\gamma\hskip .3in\text{ and }\hskip .3in H^{0}(\Qdot)_\bold 0\cong
\operatorname{coker}\gamma.
$$

\vskip .5in

We can now describe Verdier dualizing. If a perverse sheaf, $\Qdot$, on $S$ has its M-V triangle given by

\vbox{
$$ V@>\hskip .1in\operatorname{id}-\nu\hskip .1in>>V
$$
\vskip  -.2in
$$
\hbox{}\hskip -0.0in\gamma\searrow\hskip .27in\nearrow\delta
$$
\vskip -.2in
$$
\hbox{}\hskip 1in W \hskip 1in ,
$$ }
\noindent then the M-V triangle of the Verdier dual of $\Qdot$, $\Cal D\Qdot$,  is isomorphic to the dual triangle
\vskip .2in
\vbox{
$$ V^*@>\hskip .1in\operatorname{id}-\nu^t\hskip .1in>>V^*
$$
\vskip  -.2in
$$
\hbox{}\hskip -0.05in\delta^t\searrow\hskip .4in\nearrow\gamma^t
$$
\vskip -.2in
$$
\hbox{}\hskip 1in W^* \hskip 1in .
$$ }

\vskip .5in

Now, we return to the situation of the introduction, and consider the case where $\Qdot=\Pdot:=\big(\phi_f\Bbb {\frak K}^\bullet_{{}_{\Cal
U}}\big)_{|_S}[n]$. The M-V triangle of $\Pdot$,

\vskip .1in

\vbox{
$$ V@>\hskip .1in\operatorname{id}-\nu\hskip .1in>>V
$$
\vskip  -.2in
$$
\hbox{}\hskip -0.0in\gamma\searrow\hskip .27in\nearrow\delta
$$
\vskip -.2in
$$
\hbox{}\hskip 1in W \hskip 1in ,
$$ }
\noindent has $V\cong\frak K^{\lambda^1_f}$ and $W\cong \frak K^{\lambda^0_f}$, where $\lambda^0_f$ is the zero-dimensional L\^e number
(again, see [{\bf M2}]). This $0$-dimensional L\^e number can calculated effectively; let $\Gamma^1_{f, z_0}$ denote the relative polar
curve, then we have 
$\lambda^0_f= \left(\Gamma^1_{f, L}\cdot V\Big(\frac{df}{dL}\Big)\right)_\bold 0$. Note that 
$$\operatorname{ker}\gamma\ \cong\ \widetilde
H^{n-1}(\mfz) \hskip .3in\text{ and }\hskip .3in  \operatorname{coker}\gamma\ \cong\ \widetilde
H^{n}(\mfz),
$$

\noindent and the self-duality of
$\Pdot=\big(\phi_f\Bbb {\frak K}^\bullet_{{}_{\Cal U}}\big)_{|_S}[n]$ is equivalent to saying that its M-V triangle is
(non-naturally) isomorphic to its dual.

\vskip .4in

\noindent\S2. {\bf The Morse Modification}  

\vskip .2in

In this section, we shall place ourselves in a more general setting, and show how the Morse-theoretic result of L\^e in [{\bf L}] actually
yields an important method of modifying the sheaf of vanishing cycles.

\vskip .1in

Throughout this section, we continue with $\Cal U$ being an open neighborhood of the origin in $\Bbb C^{n+1}$, and let $g:(\Cal U, \bold
0)\rightarrow(\Bbb C, 0)$ be a complex analytic function (with an arbitrary critical locus).  Let $L:\Bbb C^{n+1}\rightarrow\Bbb C$ be a
linear form.

L\^e's attaching formula in [{\bf L}] provides data as to how the Milnor fibre,
$F_{g,
\bold 0}$, is built from the Milnor fibre, $F_{g_{|_{V(L)}},
\bold 0}$, of the hyperplane slice:

\vskip .3in

\noindent{\bf Theorem 2.1} (L\^e, [{\bf L}]). {\it For a generic choice of $L$, the relative cohomology $H^*(F_{g,\bold 0},
F_{g_{|_{V(L)}},
\bold 0};\ \Bbb Z)$ is zero in all degrees except, possibly, in degree $n$. In addition, $H^n(F_{g,\bold 0}, F_{g_{|_{V(L)}},
\bold 0};\ \Bbb Z)$ is free Abelian of rank equal to the intersection number of the polar curve $\Gamma^1_{g, L}$ and the hypersurface
$V(g)$ at the origin.  Furthermore, for all $\bold x\in V(g)\cap V(L)-\{\bold 0\}$ near $\bold 0$, $H^*(F_{g,\bold x}, F_{g_{|_{V(L)}},
\bold x};\ \Bbb Z)$ is zero in all degrees. }

\vskip .3in

In our own paper [{\bf M}], we generalize the above result to the case where the coefficients are an arbitrary bounded, constructible
complex of sheaves. Applying that result to the special case of perverse sheaves (e.g., the shifted constant sheaf $\frak
K^\bullet_{{}_{\Cal U}}[n+1]$), we obtain:

\vskip .3in

\noindent{\bf Theorem 2.2} (Massey, [{\bf M}]). {\it Let $\Qdot$ be a perverse sheaf of $\frak K$-vector spaces on $\Cal U$. For a generic
choice of
$L$, the relative hypercohomology $\Bbb H^*(F_{g,\bold 0}, F_{g_{|_{V(L)}},
\bold 0};\ \Qdot)$ is zero in all degrees except, possibly, in degree $0$. Furthermore, for all $\bold x\in V(g)\cap V(L)-\{\bold 0\}$ 
near $\bold 0$, $\Bbb H^*(F_{g,\bold x}, F_{g_{|_{V(L)}},
\bold x};\ \Qdot)$ is zero in all degrees.

}

\vskip .4in

\noindent{\it Remark 2.3}. In the situation of Theorem 2.2, the results of [{\bf M}] also tell one how to calculate the dimension over
$\frak K$ of $\Bbb H^n(F_{g,\bold 0}, F_{g_{|_{V(L)}}, \bold 0};\ \Qdot)$.  As in Theorem 2.1, the intersection number
\hbox{$\big(\Gamma^1_{g, L}\cdot V(g)\big)_\bold 0$} appears, but now one must calculate this intersection number for various strata,
multiply by contributions from the normal data to strata (the {\it Morse modules} of strata), and then take the sum over all strata. We
shall not need that calculation in this paper.

\vskip .4in

We wish to encode the result in Theorem 2.2 in a complex of sheaves. 

\vskip .1in

Let $\Qdot$ be a perverse sheaf of $\frak K$-vector spaces on $\Cal U$. Fix a linear form
$L$. Let
$j:\Cal U\cap V(L)\hookrightarrow\Cal U$ denote the closed inclusion, and let $i:\Cal U-\Cal U\cap V(L)\hookrightarrow\Cal U$ denote the
open inclusion. In the derived category, there is a fundamental distinguished triangle 
$$ j_*j^*\Qdot[-1]\rightarrow i_!i^!\Qdot\rightarrow\Qdot@>[1]>> j_*j^*\Qdot[-1];\tag{$\dagger$}
$$ the associated long exact sequence on hypercohomology is the relative hypercohomology long exact sequence of the pair $(\Cal U, \Cal
U\cap V(L))$:
$$
\dots\rightarrow\Bbb H^{i-1}(\Cal U\cap V(L);\ \Qdot)\rightarrow \Bbb H^{i}(\Cal U,\Cal U\cap V(L);\ \Qdot)\rightarrow \Bbb H^{i}(\Cal U;\
\Qdot)\rightarrow\Bbb H^{i}(\Cal U\cap V(L);\ \Qdot)\rightarrow\dots
$$ In the distinguished triangle $(\dagger)$, $\Qdot$ and $i_!i^!\Qdot$ are perverse, but $j_*j^*\Qdot[-1]$ need not be. However, in the
important case where
$\Qdot=\frak K^\bullet_{{}_{\Cal U}}[n+1]$, it is, in fact, true that $j_*j^*\Qdot[-1]\cong j_*\frak K^\bullet_{{}_{\Cal U\cap V(L)}}[n]$
is perverse.

\vskip .3in

There is also a canonical distinguished triangle relating the nearby and vanishing cycles:
$$ (i_!i^!\Qdot)_{|_{V(g)}}\rightarrow\psi_g[-1](i_!i^!\Qdot)\rightarrow \phi_g[-1](i_!i^!\Qdot)@>[1]>>(i_!i^!\Qdot)_{|_{V(g)}}.\tag{$*$}
$$

\vskip .1in

Let $l:V(g)\cap V(L)\hookrightarrow V(g)$ denote the inclusion. Applying the functors $\psi_g[-1]$ and $\phi_g[-1]$ to the distinguished
triangle
$(\dagger)$, comparing parts of the resulting distinguished triangles via $(*)$, and applying Theorem 2.2, we see that, for generic linear
forms
$L$, the perverse sheaf $\phi_g[-1](i_!i^!\Qdot)$ has the following four properties:

\vskip .2in

\noindent i) there is a distinguished triangle 
$$l_*\phi_{g_{|_{V(L)}}}[-1](\Qdot_{|_{V(L)}}[-1])\ \rightarrow\ \phi_g[-1](i_!i^!\Qdot)\ @>\alpha>>\ \phi_g[-1](\Qdot)\ @>[1]>>\
l_*\phi_{g_{|_{V(L)}}}[-1](\Qdot_{|_{V(L)}}[-1]);
$$

\vskip .2in

\noindent ii) if $\bold x\in V(g)-V(L)$, the morphism $\alpha$ above induces an isomorphism on stalk cohomology 
$$H^*(\phi_g[-1](i_!i^!\Qdot))_\bold x\cong H^*(\phi_g[-1](\Qdot))_\bold x;$$

\vskip .2in

\noindent iii) if $\bold x\in V(g)\cap V(L)-\{\bold 0\}$ and $\bold x$ is near $\bold 0$, $H^*(\phi_g[-1](i_!i^!\Qdot))_\bold x = 0$;

\vskip .2in

\noindent iv) $H^*(\phi_g[-1](i_!i^!\Qdot))_\bold 0$ is zero in all degrees except, possibly, degree $0$.

\vskip .3in

\noindent{\it Remark 2.4}. Note that, as $L$ is generic, ii) and iii) imply that $\operatorname{supp}\phi_g[-1](i_!i^!\Qdot) =
\operatorname{supp}\phi_g[-1]\Qdot$ near $\bold 0$. Note also that, as mentioned in Remark 2.3, it is possible to give a formula for the
dimension
${\operatorname{dim}}_{\frak K}H^0(\phi_g[-1](i_!i^!\Qdot))_\bold 0$ in terms of intersection numbers and Morse modules of strata.

\vskip .4in

\noindent{\bf Definition 2.5}. We refer to $\phi_g[-1](i_!i^!\Qdot)$ as the {\it Morse modification of $\phi_g[-1](\Qdot)$, with respect
to $L$ at
$\bold 0$}. 

We let $\overline{\Sigma_{{}_{\Qdot}}g}:=\operatorname{supp}\phi_g[-1]\Qdot$, and we refer to the restriction of the distinguished triangle in
i) to
$\overline{\Sigma_{{}_{\Qdot}}g}$ as the {\it the Morse triangle of $\phi_g[-1](\Qdot)$, with respect to $L$ at
$\bold 0$}.

\vskip .5in

Let us return now to the situation of the introduction, where $f:\Cal U\rightarrow\Bbb C$ is a complex analytic function whose critical
locus,
$\Sigma f$, equals $S=\Cal U\cap (\Bbb C\times\{\bold 0\})$, i.e., $\operatorname{supp}\phi_f\Bbb {\frak K}^\bullet_{{}_{\Cal U}} = \Cal
U\cap (\Bbb C\times\{\bold 0\})$. Let $\Qdot={\frak K}^\bullet_{{}_{\Cal U}}[n+1]$. Then, $\Qdot$ is perverse, and
$\overline{\Sigma_{{}_{\Qdot}}f}$ equals the line segment $S$. Thus, the Morse triangle of $\phi_f[-1](\Qdot)$ is a short exact sequence in
the category of perverse sheaves on $S$ and, hence, by our discussion in Section 1, is equivalent to a short exact sequence of M-V triangles:

\vskip .4in
\noindent{\bf Theorem 2.6}. {\it There is a short exact sequence of M-V triangles

\vskip .2in

\vbox{
$$
\hbox{}\hskip .1in 0@>\hskip .47in>> 0\hskip .55in \frak K^{{}^{\lambda^1_f}}@>\hskip .1in\operatorname{id}-\nu\hskip .1in>>\frak
K^{{}^{\lambda^1_f}}\hskip .54in \frak K^{{}^{\lambda^1_f}}@>\hskip .1in\operatorname{id}-\nu\hskip .1in>>\frak K^{{}^{\lambda^1_f}}
$$
\vskip  -.2in
$$ 0\hskip .05in\rightarrow\hskip .2in\hbox{}\hskip -0.0in\searrow\hskip .3in\nearrow\hskip .3in\rightarrow\hbox{}\hskip
.3in\theta\searrow\hskip .3in\nearrow\omega\hskip .3in\rightarrow\hbox{}\hskip .3in\gamma\searrow\hskip .3in\nearrow\delta\hskip
.2in\rightarrow \hskip .05in 0,
$$
\vskip -.2in
$$
\hbox{}\hskip .73in \frak K^{\mu_0}\hskip 1.25in  \frak K^{\zeta}\hskip 1.45in  \frak K^{{}^{\lambda^0_f}}\hskip .8in\hbox{}
$$ }

\vskip .1in 

\noindent where, as in the introduction and in Section 1, $\nu$ is the internal monodromy, $\lambda^1_{f}$ is the $1$-dimensional L\^e
number, 
$\lambda^0_f$ is the $0$-dimensional L\^e number, $\mu_0$ is the Milnor number of the isolated critical point (at
the origin) of $f_{|_{V(L)}}$, $\zeta:= \mu_0 + \lambda^0_f$, and $\theta$ is an injection.}

\vskip .3in

 The fact that $\theta$ is an {\bf injection} follows from property iv) of the Morse modification.
We let
$\beta$ denote the linear map from $\frak K^{\mu_0}$ to  $\frak K^{\zeta}$.

\vskip .4in

\noindent{\bf Corollary 2.7}. {\it $\operatorname{dim}_{\frak K}\widetilde H^{n-1}(\mfz)\leqslant\lambda^1_f$, with equality holding if
and only if $\operatorname{im}\theta\subseteq \operatorname{im}\beta$. }

\vskip .2in

\noindent{\it Proof}. As we mentioned in Section 1, $\widetilde H^{n-1}(\mfz)\cong \operatorname{ker}\gamma$; thus, certainly, the
inequality holds. In addition, it follows that the equality holds if and only if $\gamma$ is the zero map. As $\theta$ is an injection,
$\gamma=0$  is equivalent to $\operatorname{im}\theta\subseteq \operatorname{im}\beta$.
\qed

\vskip .4in

\noindent\S3. {\bf The Main Theorem}  

\vskip .2in

We return to the setting from the introduction, where we are given a complex analytic function $f:(\Cal U, \bold 0)\rightarrow(\Bbb C,
0)$, whose critical locus at the origin is a line. We may select a generic linear form $L$, and consider the one-parameter family of
isolated singularities given by
$f_t:=f_{|_{V(L-L(t))}}$, and also consider their corresponding Milnor numbers $\mu_t$. Viewing the situation in this from, we see that
this $\mu_0$ agrees with the $\mu_0$ of Theorem 2.6, and that
$\lambda^1_f=\mu_t$, for all small $t\neq 0$.

\vskip .3in

The Milnor monodromy of $f$, induced by letting the value of $f$ travel around a small circle centered at the origin, yields an
automorphism, $T$, of the entire short exact sequence of M-V triangles in \hbox{Theorem 2.6}; that is, we obtain automorphisms of each
vertex of each M-V triangle which commute with all the other maps. The main theorem will follow immediately from the existence of these
automorphisms, combined with Corollary 2.7 and the following theorem of A'Campo.

\vskip .4in

\noindent{\bf Theorem 3.1} (A'Campo, [{\bf A'C}]). {\it Let $X$ be a complex analytic space, $g:X\rightarrow\Bbb C$ a complex analytic
function, and $\Fdot$ a bounded, constructible complex of sheaves of $\frak K$-vector spaces on $X$. Let $\frak m_{{}_{X, \bold p}}$
denote the maximal ideal of $X$ at a point $\bold p\in V(g)$. If $g\in \frak m^2_{{}_{X, \bold p}}$, then the Lefschetz number of the
Milnor monodromy automorphism on the stalk cohomology of
$\psi_g(\Fdot)$ at $\bold p$ equals $0$. }

\vskip .4in

\noindent{\it Remark 3.2}. In fact, in [{\bf A'C}], A'Campo states the above result only in the case where $\frak K=\Bbb C$. However, one
sees easily that his proof is valid for arbitrary fields; by using resolution of singularities (and filtrations and long exact sequences),
one is reduced to the case where $g$ is of the form $z_0^{\alpha_0} z_1^{\alpha_1}\dots z_n^{\alpha_n}$, where the theorem certainly holds
over arbitrary fields.

\vskip .4in

Note that A'Campo's result is a statement concerning $\psi_g$, not $\phi_g$. We need to discuss the implications for the monodromy
automorphism of our short exact sequence of M-V triangles in Theorem 2.6.

The monodromy automorphisms $T_{\mu_0}:\frak K^{\mu_0}\rightarrow\frak K^{\mu_0}$ and $T_{\lambda^1_f}:\frak
K^{\lambda^1_f}\rightarrow\frak K^{\lambda^1_f}$ are, respectively, the Milnor automorphisms induced on $\widetilde H^{n-1}(F_{f_0, \bold
0};\ \frak K)$ and $\widetilde H^{n-1}(F_{f_t, \bold 0};\ \frak K)$ for small $t\neq 0$. As each $f_t$ is reduced (in a neighborhood of
$\bold 0$ in $\Cal U$) and has an isolated critical point at the origin, A'Campo's result implies that the traces are given by
$\operatorname{tr}(T_{\mu_0}) =
\operatorname{tr}(T_{\lambda^1_f}) = (-1)^n$.

\vskip .4in

With all of our preliminary discussion and results, the proof of the main theorem is now simple.

\vskip .4in

\noindent{\bf Theorem 3.3}. {\it If $\mu_0=1+\lambda^1_f$, then there are strict inequalities $$\operatorname{dim}_{\frak
K}\widetilde H^{n-1}(\mfz)<\lambda^1_f\hskip .2in\text{ and }\hskip .2in\operatorname{dim}_{\frak K}\widetilde
H^{n}(\mfz)<\lambda^0_f.$$ }

\vskip .2in

\noindent{\it Proof}. Recall that in the short exact sequence of M-V triangles in Theorem 2.6,  $\widetilde
H^{n-1}(\mfz)\cong\operatorname{ker}\gamma$ and $\widetilde H^{n}(\mfz)\cong\operatorname{coker}\gamma$. Using the result and notation
from Corollary 2.7, we see that the inequalities of the theorem hold, unless
$\operatorname{im}\theta\subseteq\operatorname{im}\beta$. 

If $\mu_0=1+\lambda^1_f$, then
$\operatorname{im}\theta$ and $\operatorname{im}\beta$ are invariant subspaces of the monodromy, whose dimension differs by one, such that
the trace of the monodromy on each of the subspaces is $(-1)^n$. If $\operatorname{im}\theta\subseteq\operatorname{im}\beta$, then the
induced monodromy isomorphism on $\operatorname{im}\beta/\operatorname{im}\theta\cong\frak K$ would have trace equal to zero; this is
impossible.\qed

\vskip .4in

\noindent\S4. {\bf Concluding Remarks}  

\vskip .2in

Theorem 3.3 may leave the reader asking several questions:

\vskip .1in

\noindent $\bullet$ Why care about a result which such a restrictive hypothesis?

\vskip .1in 

\noindent $\bullet$ Are all of the abstract tools of this paper really necessary to prove Theorem 3.3?

\vskip .1in

\noindent $\bullet$ Does Theorem 3.3 yield the final, general restriction on the cohomology of the Milnor fibre of an affine line
singularity?

\vskip .1in

We address these questions below.

\vskip .4in

\noindent {\bf Why care about a result which such a restrictive hypothesis?}

\vskip .1in

It is notoriously difficult to prove any general statement concerning the cohomology of Milnor fibres of hypersurfaces with non-isolated
singularities. It is especially difficult to prove such results where the hypotheses are so easy to calculate effectively, as is the
hypothesis of Theorem 3.3. 

For example, recall from the introduction that there is an injection
$$
\widetilde H^{n-1}(\mfz)\hookrightarrow \operatorname{ker}(\operatorname{id}-\nu).
$$ Thus, if the internal monodromy, $\nu$, is not the identity, then we could again conclude the result of Theorem 3.3, that
$$\operatorname{dim}_{\frak K}\widetilde H^{n-1}(\mfz)<\lambda^1_f.$$

However, it is not so easy, in general, to decide when $\nu$ is the identity. Moreover, it is {\bf not} true that Theorem 3.3 results from
forcing the monodromy to be non-trivial, i.e., it is {\bf not} true that: if
$\mu_0=1+\lambda^1_f$, then $\nu$ is not the identity.  Consider $f_t:= y^2-x^3-t^2x^2$; the reader may verify that the internal
monodromy is the identity,
$\lambda^1_f=1$, $\mu_0=2$, $\lambda^0_f=5$, and the Milnor fibre at the origin has the homotopy-type of a bouquet of four
$2$-spheres.

\vskip .4in

\noindent{\bf Are all of the abstract tools of this paper really necessary to prove Theorem 3.3?}

\vskip .1in

This is more difficult to answer. Our proof uses complexes of sheaves to compare the Milnor fibres of $f$, $f_0$, and $f_t$ (for $t\neq
0$) and their respective monodromies. It may be that there is a more elementary proof. If there is, we do not know of one. Moreover, our
proof is really incredibly simple, once one understands the background material. In addition, we believe that the Morse modification and
the Morse triangle of Definition 2.5 will prove useful in studying other questions concerning non-isolated hypersurface singularities.

\vskip .4in

\noindent{\bf Does Theorem 3.3 yield the final, general restriction on the cohomology of the Milnor fibre of an affine line singularity?}

\vskip .1in

Unfortunately, the answer is: no. This paper was motivated by a result of Siersma in [{\bf S}]. Translating Sierma's result into the
language and notation of this paper, it says: if $\lambda^1_f = 1$ and $\mu_0\neq 1$, then $\operatorname{dim}_{\frak K}\widetilde
H^{n-1}(\mfz)<\lambda^1_f$, i.e.,
$\widetilde H^{n-1}(\mfz)=0$. In terms of M-V triangles, Siersma's result says that, if $\lambda^1_f = 1$ and $\mu_0\neq 1$, then
the M-V triangle of $\big(\phi_f\Bbb {\frak K}^\bullet_{{}_{\Cal U}}\big)_{|_S}[n]$ is {\bf not} the direct sum 

\vbox{
$$ 0@>\hskip .3in>> 0\hskip .55in \frak K@>\hskip .2in\hskip .1in>>\frak K
$$
\vskip  -.2in
$$
\hbox{}\hskip -.07in\searrow\hskip .16in\nearrow\hskip .27in\oplus\hskip .27in\searrow\hskip .2in\nearrow
$$
\vskip -.2in
$$
\hbox{}\hskip .07in\frak K^{\lambda^0_f}\hskip 1in 0\hskip .1in\hbox{}.
$$ }

\vskip .2in

Our result does not imply that of Siersma or, if it does, we do not see how. However, our result and Siersma's have a similar feel, and it
is easy to believe that there is one statement and proof that yields both results as special cases of a more general result.

On the other hand, it may be that Theorem 3.3 can, in fact, be used to recover the result of Siersma. The assumption that
$\mu_0=1+\lambda^1_f$ is equivalent to the assumption that the generic reduced relative polar curve, $\Gamma^1_{f, L}$, has a
single smooth component through the origin. In the general case where the polar curve is itself singular, it may be possible to
perform an embedded resolution of $\Gamma^1_{f, L}$, or of 
$\Sigma\cup \Gamma^1_{f, L}$, and then use Theorem 3.3 to gain information about the Milnor fibre of $f$. Thus far, we have not succeeded
with this approach.

\enddocument